\documentclass[10pt,final,a4paper]{article}
\usepackage[colorinlistoftodos, textwidth=2cm, shadow]{todonotes}
\usepackage{amssymb}
\usepackage{amsfonts}
\usepackage{gensymb}
\usepackage{threeparttable}
\usepackage{amsthm}
\usepackage[cmex10]{amsmath}
\usepackage{hyperref}
\usepackage[affil-it]{authblk}

\begin{document}

\title{Aggregated motion estimation for real-time MRI reconstruction}

\author{Housen~Li$^{1,2}$\thanks{Corresponding author: housen.li@mpibpc.mpg.de},
        Markus~Haltmeier$^{3}$,
        Shuo~Zhang$^{4}$,
        Jens~Frahm$^{4}$,
        and~Axel~Munk$^{1,2}$,%
        }\affil{$^1$Institute for Mathematical Stochastics, University of G\"ottingen, G\"ottingen, Germany. \\
$^2$Max Planck Institute for Biophysical Chemistry, G\"ottingen, Germany. \\
$^3$Department of Mathematics, University of Innsbruck, Innsbruck, Austria. \\
$^4$Biomedizinische NMR Forschungs GmbH, Max Planck Institute for Biophysical Chemistry, G\"ottingen, Germany.}

\date{}
\maketitle

\begin{abstract} 
Real-time magnetic resonance imaging (MRI) methods generally shorten the measuring time by acquiring less data than needed according to the sampling theorem. In order to obtain a proper image from such undersampled data, the reconstruction is commonly defined as the solution of an inverse problem, which is regularized by \textit{a priori} assumptions about the object. While practical realizations have hitherto been surprisingly successful, strong assumptions about the continuity of image features may affect the temporal fidelity of the estimated images. Here we propose a novel approach for the reconstruction of serial real-time MRI data which integrates the deformations between nearby frames into the data consistency term. The method is not required to be affine or rigid and does not need additional measurements. Moreover, it handles multi-channel MRI data by simultaneously determining the image and its coil sensitivity profiles in a nonlinear formulation which also adapts to non-Cartesian (e.g., radial) sampling schemes. Experimental results of a motion phantom with controlled speed and \textit{in vivo} measurements of rapid tongue movements demonstrate image improvements in preserving temporal fidelity and removing residual artifacts.
\end{abstract}


\section{Introduction}

Imaging speed is crucial in real-time MRI studies of physiologic processes, ranging from cardiovascular imaging to noninvasive monitoring of interventional (surgical) procedures. Because the physical acceleration of the data acquisition process is limited by physiologic regulations to prevent peripheral nerve stimulation, most strategies ultimately reduce the measuring time by acquiring less data, while attempting to preserve the quality of the reconstructed image. A first development along this line is the \emph{parallel imaging} concept which takes advantage of multiple receiver coils to acquire data simultaneously. Such techniques benefit from encoding part of the spatial information of an object into spatially complementary coil sensitivities, which are generally unknown and also depend on the actual experimental condition. Therefore, coil sensitivity profiles are either explicitly pre-calibrated in image space, like SENSE \cite{Pruessmann1999}, or implicitly in $k$-space, like SMASH \cite{Sodickson1997} and GRAPPA \cite{Griswold2002}. Unfortunately, however, such pre-calibration techniques make only suboptimal use of the available data from multiple receiver channels, so that respective errors in the estimated coil profiles may lead to artifacts in the iteratively optimized image -- already for moderate acceleration factors of about two to three. An improved strategy is to compute spin density maps and coil profiles at the same time by means of a nonlinear formulation of the inverse reconstruction problem \cite{Ying2007,Uecker2008}. In this case, when aiming for high temporal resolution, the use of strongly undersampled data introduces additional ill-posedness to the reconstruction problem. In order to stabilize the problem and obtain plausible solutions, it is necessary to incorporate  \textit{a priori} information about the unknown object (described by its spin density) and the coil profiles into the reconstruction method. In \cite{Uecker2010,Uecker2010a,Zhang2010a}, a temporal $L^2$-regularization on the object was studied, which is, however, usually too weak to remove residual artifacts. For example, temporally flickering artifacts are observed for a radial sampling scheme which employs complementary sets of spatial encodings in consecutively acquired datasets. In practice, a temporal median filter may effectively diminish residual streaking artifacts though at the expense of degrading the true temporal resolution. Alternatively, the total variance or total generalized variance were used for regularization \cite{Knoll2012}, which also reduce streaking artifacts, but fail to recover small-scale details of the object and therefore sacrifice spatial resolution. In general, regularization methods alone seem to be unable to provide artifact-free images with high spatial and temporal resolution in real-time MRI scenarios with pronounced undersampling.

An alternative strategy for improved image quality is to integrate information about moving object features into the reconstruction by exploiting multiple measurements at different (neighboring) time points. One of the most effective means for motion compensation in MRI is the `navigator echo' technique and its variants \cite{Ehman1989,Ward2000,Welch2001,Kouwe2006,White2009,Ooi2009,Lin2009,Nielsen2011}. In these methods, a navigator signal is repetitively acquired during the scan to extract specific motion information. The need to insert multiple navigator modules into the MRI sequence may be avoided by `self-navigating' techniques which may determine motions from the actual data. Early applications \cite{Schaffter1999,Batchelor2005} used partial or full $k$-space data in a block-based or parametric manner, but failed to detect complex motions such as elastic deformations. Lately, more flexible motion-detection techniques were developed for free-breathing cine MRI studies \cite{Kellman2008,Kellman2009,Hansen2012,Usman2012}, but they can only compensate for slow (e.g., respiratory) movements that affect a faster (e.g., cardiac) motion of interest. Motion compensation was also combined with conventional parallel MRI reconstructions \cite{Odille2008,Wei2012}.

To overcome the aforementioned limitations, this work presents a novel reconstruction method for real-time MRI by integrating the idea of a `self-navigating' motion into a nonlinear formulation of the inverse problem which simultaneously estimates spin density and coil sensitivity profiles. Based on a non-parametric motion estimation, the new method generates images with high temporal fidelity and reduced residual artifacts. The validation of the approach in comparison to current algorithms employed a motion phantom with controlled speed as well as real-time MRI studies of movements during human sound production.

\section{Theory}
The proposed method is based on a recently developed reconstruction from highly undersampled radial MRI acquisitions with multiple receiver coils \cite{Uecker2008,Uecker2010}. It generalizes the respective data consistency term to incorporate an aggregated reconstruction from multiple frames with non-parametric motion correction (i.e., AME = \emph{Aggregated Motion Estimation}) and is schematically outlined in Figure \ref{workflow}.

\subsection{Real-time Magnetic Resonance Imaging}
The real-time MRI data acquisition of the $t$-th ($t \in \mathbb{N}$)  frame from multiple receiver coils is given by
\begin{equation}\label{meas_equ}
y_{t,l} = S_t \mathcal{F}(\rho_t \cdot c_{t,l}) + \varepsilon_{t,l}, \quad l \in \Lambda (:= \{1, \ldots, N\}),
\end{equation}
where $\rho_t$ denotes the spin density, $c_{t,l}$ the sensitivity profiles of $N$ individual receiver coils, $\varepsilon_{t,l}$  the noise, $S_t: f \mapsto (f(k_{t,j}))_{j=1}^{M} $  the sampling operator at the positions $(k_{t,j})_{j=1}^{M}$ in $k$-space, and $\mathcal{F}: f\mapsto \int f(x) e^{-2\pi i \langle x, k \rangle} dx$  the Fourier transform. The goal is to obtain a serial stream of images $(\rho_t)_{t\in \mathbb{N}}$ with high spatial and temporal resolution from the measured data $(y_{t,l})_{t\in N, l\in \Lambda}$. The MRI measuring time per frame is mainly determined by the number of samples $M$ (times the physical repetition time TR needed for radiofrequency excitation and spatial encoding), which therefore is kept as small as possible. On the other hand, for pronounced undersampling conditions, equation (\ref{meas_equ}) becomes increasingly ill-conditioned. As a consequence, the inversion of the system leads to an amplification of noise which in turn results in low-resolution images. Thus, a proper choice of $M$ should be a sensible trade-off between temporal and spatial resolution.

\subsection{Aggregated Motion Estimation for Nonlinear Reconstruction}
For explanation purpose, in this subsection, we first assume that the motion (or deformation) $\phi_{t,s}$ from $\rho_t$ to $\rho_{t+s}$, i.e. $\phi_{t,s}(\rho_t)\approx\rho_{t+s}$, is known for every $s \in K$, with some $K \subset \mathbb{Z}$. For example, $K$ can be $\{-2, -1, 0, 1, 2\}$, which corresponds to 5 successive acquisitions. In this case, only two future frames are included corresponding to an ignorable waiting time of around 60 ms in our applications.  By variable substitution, the spin density of the $t$-th frame satisfies

\begin{equation}\label{ame_equ}
S_{t+s} \mathcal{F}(\phi_{t,s}(\rho_t) \cdot c_{t+s,l}) = y_{t+s,l}, \quad l \in \Lambda, s\in K.
\end{equation}

Thus, if successive frames rely on complementary data samples in $k$-space, the reconstruction takes advantage of $|K|\cdot M$ samples for recovering $\rho_t$, while the temporal resolution remains unchanged, i.e. corresponding to  $M$ times the repetition time TR. Accordingly, while keeping the temporal resolution, the approach may obtain images with higher spatial resolution from (\ref{ame_equ}). For sake of clarity, we rewrite (\ref{ame_equ}) by an abstract nonlinear operator equation

\begin{equation*}
F_{t+s}(\Phi_{t,s}(x)) = y_{t+s}, \mbox{ for } s \in K, \mbox{ with } x = \left(\begin{array}{c}
\rho_t\\
(c_{t+s',l})_{s'\in K, l \in \Lambda}
\end{array}\right), \  y_{t+s} = (y_{t+s,l})_{l \in \Lambda}.
\end{equation*}
Here, $F_{t+s}: x \mapsto \left(S_{t+s}\mathcal{F}(\rho_t \cdot c_{t+s,l})\right)_{l \in \Lambda}$ and $\Phi_{t,s}: x \mapsto \left(\phi_{t,s}(\rho_t), (c_{t+s',l})_{s' \in K, l \in \Lambda}\right)^{T}$, for every $s \in K$.

These equations are solved for the unknown $x$ by a Newton-type method, whose key idea consists in repeatedly linearizing the operator equation $F_{t+s}(\Phi_{t,s}(x)) = y_{t+s}$, $s \in K$, around some approximate solution $x_n$, and solving the linearized problems
\begin{equation}\label{lin_equ}
F_{t+s}'(\Phi_{t,s}(x_n))\Phi_{t,s}'(x_n)(x - x_n) = y_{t+s} - F_{t+s}(\Phi_{t,s}(x_n)), \quad s \in K, 
\end{equation}
As the real-time MRI problem is highly ill-posed, $\left(F_{t+s}'(\Phi_{t,s}(x_n))\Phi_{t,s}'(x_n)\right)^{-1}$ is not bounded or seriously ill-conditioned. The standard Newton method is not applicable and may not even be well defined for noise-free data, because $y_{t+s}-F_{t+s}(\Phi_{t,s}(x_n))$ is not guaranteed to lie in the range of $F_{t+s}'(\Phi_{t,s}(x_n))\Phi_{t,s}'(x_n)$. Therefore, some regularization method has to be employed  for solving the linearized equation (\ref{lin_equ}).

Because only the product of the spin density and coil profiles is determined, the real-time MRI problem is undetermined even in the fully sampled case. Although the image may contain fine structures, the coil profiles are generally rather smooth. As in \cite{Uecker2008}, this can be ensured by introducing a term promoting smoothness which may be given by a Sobolev norm $\|f\|_m := \|(1+b\|\cdot\|^2)^{m/2}(\mathcal{F}f)(\cdot)\|$, $m\in \mathbb{R}_{+}$. It penalizes high spatial frequencies by a polynomial of degree $m$ as a function of the distance to the centre of $k$-space. The object (i.e., its spin density) usually deforms continuously and smoothly from frame to frame, so an efficient regularization penalizes the differences between neighboring frames to ensure temporal continuity. By combining the standard Newton method and the aforementioned regularization, we obtain the well-known \emph{iteratively regularized Gauss-Newton method} (IRGNM) \cite{Bakushinskii1992,Bauer2009} for solving (\ref{ame_equ})
\begin{eqnarray*}
h_n & = & \mathop{\arg\min}\limits_{h}\{ \sum_{s \in K}\|F_{t+s}'(\Phi_{t,s}(x_n))\Phi_{t,s}'(x_n)h - \left( y_{t+s} - F_{t+s}(\Phi_{t,s}(x_n))\right)\|^2 \\
&&{+}\: \alpha_n\sum_{s\in K, l \in \Lambda}\|a(1+b\|\cdot\|^2)^\frac{m}{2}\mathcal{F}(c_{t+s,l}^{(n)}+h_{c_{t+s,l}}-c^0_{t+s,l})\|^2,\\
&&{+}\: \alpha_n\|\rho_t^{(n)}+h_{\rho_t}-\rho_t^0\|^2\} \\
x_{n+1} & = & x_n + h_n,
\end{eqnarray*}
where $x_n = (\rho_t^{(n)}, (c_{t+s,l}^{(n)})_{s\in K, l \in \Lambda})^T$, $h = (h_{\rho_t}, (h_{c_{t+s,l}})_{s \in K, l \in \Lambda})^T$ and $x_0 = (\rho_t^0, (c_{t+s,l}^{0})_{s\in K, l \in \Lambda})^T$ the initial guess. If the initialization is close enough to the true solution, the choice of $\alpha_n := \alpha_0 q^n$, with $0 < q < 1$, is usually sufficient for convergence (cf. \cite{Hohage1999}). Because the IRGNM method reduces to a Gauss-Newton method for $\alpha_n = 0$, it uses the more robust descent direction at the beginning of the iterative process (far from the solution) and the faster convergent algorithm at the end (near the solution). The choice of parameter $a$ serves as a balance between penalization of the spin density and coil profiles. The quadratic optimization is equivalent to solve its normal equation, precisely a linear equation with a symmetric coefficient matrix. Unfortunately, this linear equation is numerically ill-conditioned for large $m$. A simple preconditioning by the following variable substitution can significantly reduce its condition number, making it numerically stable. Let
\[
M: \left(\begin{array}{c}
\rho_{t}\\
(c_{t+s,l})_{s\in K, l \in \Lambda}
\end{array}\right) \mapsto
\left(\begin{array}{c}
\rho_t\\
(Wc_{t+s,l})_{s\in K, l \in \Lambda}
\end{array}\right),\]
with $W = \mathcal{F}^*\circ a^{-1}(1+b\|\cdot\|^2)^{-m/2}$. If $\tilde{x} := M^{-1}x,\, G_{t+s} := F_{t+s}\circ\Phi_{t,s}\circ M,\, s \in K$, then an equivalent form of IRGNM is given by
\begin{eqnarray*}
\tilde{h}_n & = & \mathop{\arg\min}\limits_{\tilde{h}} \sum_{s \in K}\|G_{t+s}'(\tilde{x}_n)\tilde{h} - (y_{t+s} - G_{t+s}(\tilde{x}_n))\|^2 + \alpha_n \|\tilde{x}_n+\tilde{h}-\tilde{x}_0\|^2\\
\tilde{x}_{n+1} & = & \tilde{x}_n + \tilde{h}_n.
\end{eqnarray*}
Explicitly, the optimality condition for this quadratic optimization is
\begin{equation}\label{final_NE}
\left(\sum_{s \in K}G_{t+s}'(\tilde{x}_n)^*G_{t+s}'(\tilde{x}_n)+\alpha_n I\right)\tilde{h}_n = \sum_{s \in K}G_{t+s}'(\tilde{x}_n)^*\left(y_{t+s} - G_{t+s}(\tilde{x}_n)\right) - \alpha_n(\tilde{x}_n - \tilde{x}_0).
\end{equation}

\begin{figure}[!t]
\centering
\includegraphics[width=0.8\textwidth]{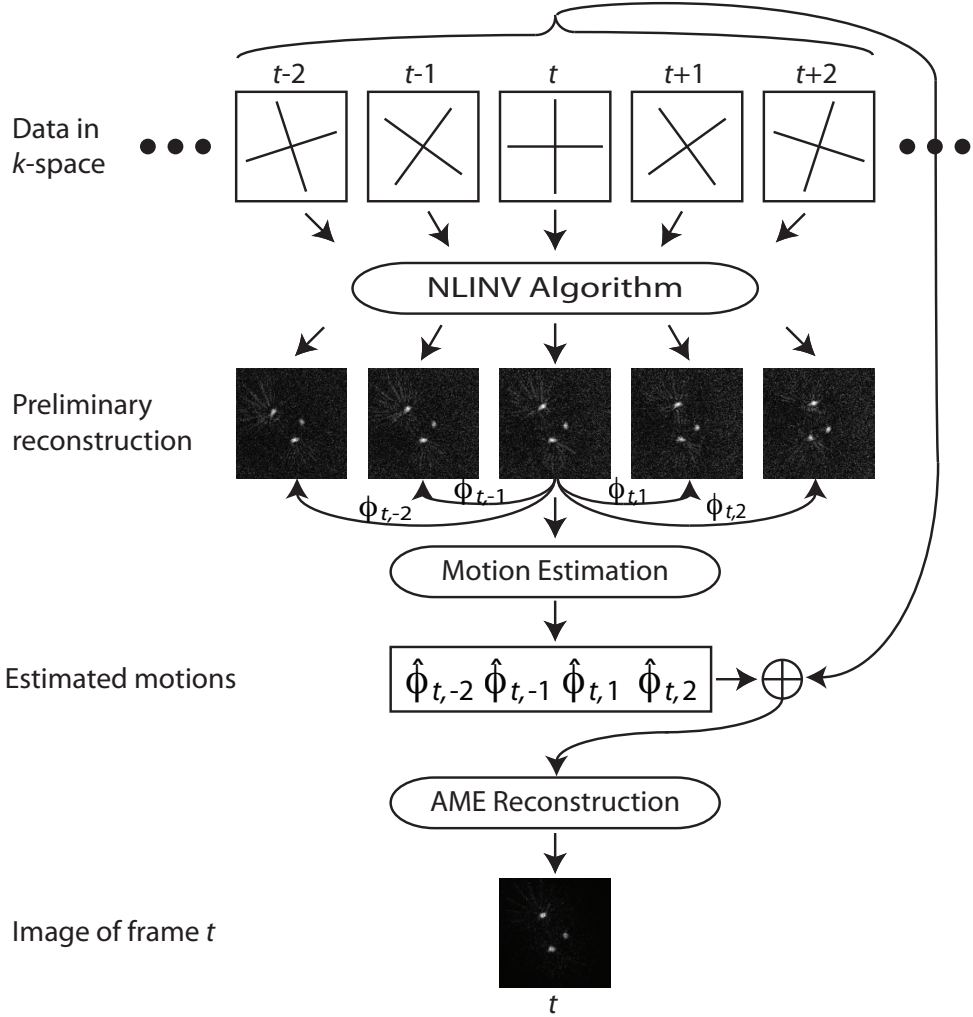}
\caption{Flow diagram illustrating the AME-based nonlinear reconstruction method for reconstructing the $t$-th frame, with $K = \{-2, -1, 0, 1, 2\}$.}
\label{workflow}
\end{figure}

Its discretized form can efficiently be solved by the \emph{conjugate gradient} (CG) algorithm. Together with the motion estimation described in the next subsection, this strategy represents our novel AME-based nonlinear reconstruction method for high temporal and spatial resolution (see Figure  \ref{workflow}). In this manner, multiple acquisition frames are exploited for reconstruction with proper motion correction, implicitly increasing the sampling rate while preserving temporal sharpness. 

\subsection{Motion Estimation}
In general, the differential positional displacements between nearby frames are not known. Therefore, we first precompute each frame by the nonlinear inversion (NLINV) method introduced in \cite{Uecker2008,Uecker2010}, which corresponds to $K = \{0\}$ as defined in the previous subsection. Subsequently, the motion is estimated on these precomputed images, using the TV-$L^1$ optical flow model (cf. \cite{Wedel2009,Chambolle2011}). In detail, the motion $\phi_{t,s}$ from $\rho_t$ to $\rho_{t+s}$ is estimated by $\hat{\phi}_{t,s}(\rho_t)(x) = \rho_t(x+\hat{u}(x))$, with $\hat{u}$ given by the solution to
\begin{equation*} 
\mathop{\min}\limits_{u, v}  \|\rho_{t} + \nabla \rho_t \cdot u- \rho_{t+s} + v \|_1 + \lambda \|\nabla u\|_1 + \mu \|\nabla v\|_1.
\end{equation*}

Because $\phi_{t,s}(\rho_t)(x) := \rho_t(x+u(x)) \approx \rho_t + \nabla \rho_t \cdot u$, the auxiliary variable $v$ models the varying reconstruction artifacts in different frames. For radial MRI acquisitions, the residual streaking artifacts have a relatively low total variation in comparison to the object which contains all local structures. Therefore, it is expected that $u$ can only capture the true motion of the object instead of the undesirable motion of the artifacts contained in the precomputed images. In order to avoid the impact of outliers on motion estimation, we only used the $L^1$ norm. This non-smooth minimization can efficiently be solved by the first-order primal-dual algorithm proposed in \cite{Chambolle2011}.

\subsection{Discretization} \label{discrete}
By denoting $C_{t+s,l}(\tilde{x}):=\phi_{t,s}(\rho)\cdot W\tilde{c}_{t+s,l}$, with $\tilde{x} := (\rho,(\tilde{c}_{t+s',l'})_{s'\in K,l' \in\Lambda})^T$, for $s\in K$, $l\in \Lambda$, a detailed formula of (\ref{final_NE}) is
\begin{eqnarray*}
&&\left(\sum_{s\in K, l \in \Lambda} C_{t+s,l}'(\tilde{x}_n)^*\mathcal{F}^*S_{t+s}^*S_{t+s}\mathcal{F}C_{t+s,l}'(\tilde{x}_n) + \alpha_n I\right)\tilde{d}_n\\
&=&\sum_{s\in K, l \in \Lambda} C_{t+s,l}'(\tilde{x}_n)^*\mathcal{F}^*S_{t+s}^*\left(y_{t+s,l}-S_{t+s}\mathcal{F}C_{t+s,l}(\tilde{x}_n)\right)-\alpha_n(\tilde{x}_n-\tilde{x}_0).
\end{eqnarray*}
This equation will be solved by the CG algorithm which requires repeated application of the operations $S_{t+s}\mathcal{F}$ and $\mathcal{F}^* S_{t+s}^*$. ÊFor numerical computation, every function involved needs to be approximated by a discretized form of points on a rectangular grid. Since density and coil profiles are compactly supported, the Fourier transform $\mathcal{F}$ can be computed by fast Fourier transform (FFT) with proper periodic extension.

If the sampling trajectory represents a non-Cartesian radial scheme as used for real-time MRI in \cite{Uecker2008,Uecker2010}, the computation involving the sampling operators $S_{t+s}$ and $S_{t+s}^*$ is not straightforward. With $\langle S_{t+s}f, g \rangle = \langle f, S_{t+s}^* g \rangle$ and $S_{t+s}f = \sum_{j=1}^M f(k_{t+s,j})$, we have $S_{t+s}^*g = \sum_{j=1}^M \delta(\cdot - k_{t+s,j}) g(k_{t+s,j})$. Then, $\mathcal{F}^*S_{t+s}^*y_{t+s,l} = \sum_{j=1}^My_{t+s,l}e^{2\pi i \langle k_{t+s,j}, \cdot \rangle}$ can be computed by inverse FFT after gridding \cite{Osullivan1985,Jackson1991,Beatty2005}, or nonuniform FFT \cite{Fessler2003,Keiner2009}. With respect to $\mathcal{F}^*S_{t+s}^*S_{t+s}\mathcal{F}$, we have
\begin{eqnarray*}
(\mathcal{F}^*S_{t+s}^*S_{t+s}\mathcal{F}f)(x) &=& \int \sum_{j=1}^M\delta(y-k_{t+s,j}) e^{2\pi i \langle x, y \rangle}\int f(z) e^{-2 \pi i \langle k_{t+s,j}, z \rangle}dz dy\\
&=&\int \sum_{j=1}^M e^{2 \pi i \langle k_{t+s,j}, x - z \rangle} f(z) dz\\
&=& q * f\\
&=&\mathcal{F}^*\left((\mathcal{F}q)(\mathcal{F}f)\right),
\end{eqnarray*}
where $q(x) := \sum_{j=1}^M e^{2\pi i \langle k_{t+s,j} x\rangle}$. It can be computed by two FFTs and one inverse FFT, with $\mathcal{F}q$ given by the gridding algorithm. In the Cartesian case $\mathcal{F}q$ equals ones at measured points and zeros elsewhere. To sum up, equation (\ref{final_NE}) can numerically be solved in an efficient way.

\section{Methods}
\subsection{Data Acquisition}
The proposed reconstruction technique was evaluated for real-time MRI measurements of a motion phantom as well as for different parts of the human body \textit{in vivo}. All studies were conducted on a 3T MRI system (Siemens Magnetom TIM Trio, Erlangen, Germany). Continuous data acquisition was achieved by using a radiofrequency-spoiled radial FLASH (fast low angle shot) pulse sequence developed for real-time MRI \cite{Zhang2010}. T1-weighted images were generated by a short repetition time TR (approximately 2 ms) and a low flip angle of the RF excitation pulse (5 to 10 degree). A highly undersampled radial $k$-space encoding scheme was employed with an interleaved arrangement of spokes for five successive datasets (i.e., frames). Each single turn corresponded to a full image and contained only a small number of spokes that were equally distributed over a full 360\degree circle in order to homogeneously sample $k$-space. To prevent aliasing effects from object structures outside the selected field-of-view, a readout oversampling by a factor of two was used during data acquisition without compromising imaging speed or signal, for details see \cite{Zhang2010}. For human studies, healthy subjects with no known illness were recruited among the university students and written informed consent was obtained in all cases prior to each examination.

The motion phantom consisted of a polyacetal disc rotating with respect to its geometric center. Three water-filled tubes with approximately 10 mm diameter were fixed on the disc with a distance to the center of 25 mm, 37.75 mm, and 55 mm, respectively. The MRI signals were acquired using a 32-channel head coil (Siemens Healthcare, Erlangen, Germany) and the measurements were performed with three different rotational speeds at angular velocities of 0.5 Hz, 1.0 Hz, and 1.5 Hz, respectively.

Real-time MRI of the human body was performed in a supine position for studies of the heart and movements of the tongue during playing a plastic mouthpiece of a brass instrument. In the latter case subjects were asked to perform rapid tongue movements (\emph{staccato}) at a rate of about 5 Hz. A mid-sagittal image was chosen to cover the oropharyngolaryngeal area, while MRI signals were acquired by combining a 4-channel small flexible receiver coil (Siemens Healthcare, Erlangen, Germany) and a bilateral $2\times4$ array coil (NORAS MRI products, Hoechberg, Germany), using the same setup as previously reported for real-time MRI of speech generation \cite{Niebergall2012}. Cardiac MRI was performed during free breathing and without synchronization to the electrocardiogram \cite{Uecker2010a,Zhang2010a} using a 32-channel body coil consisting of an anterior and a posterior array with 16 elements each. Online image control employed conventional NLINV reconstructions with a post-processing temporal median filter (NLINV-MED). Details of the imaging parameters are summarized in Table \ref{data_par}.

\begin{table}
\centering
\caption{Acquisition parameters for real-time MRI.}\label{data_par}
\begin{tabular}{l|ccc}
  \hline
    Acquisition Parameter   & Motion Phantom & Tongue Movement  \\
  \hline
  \hline
  Repetition time (ms)  &2.28 & 2.2 \\
  Echo time (ms) & 1.48 & 1.4  \\
  Flip angle (\degree) & 8 & 5  \\
  Field-of-view ($\mbox{mm}^2$) & $256 \times 256$ & $192 \times 192$  \\
  Base resolution     & $128 \times 128$ & $128 \times 128$  \\
  Section thickness (mm)    & 5.0 & 10.0 \\
  Voxel size ($\mbox{mm}^3$)    & $2.0 \times 2.0 \times 5.0$ & $1.5 \times 1.5 \times 10.0$ \\
  Spokes per frame       & 9   &  15   \\
  \hline
\end{tabular}
\end{table}

\subsection{Image Reconstruction}
All reconstructions were done offline using an in-house software package written in Matlab (R2012a, The MathWorks, Natick, MA). In the first step, data from up to 32 receiver channels were combined into a small set of 10 virtual channels based on a principal component analysis, as previously described in \cite{Buehrer2007,Zhang2010a}. For the interpolation in $k$-space from radial spokes to Cartesian grids, a Kaiser-Bessel window function with $L = 6$, $\beta = 13.8551$ and a 1.5 fold oversampling was used \cite{Beatty2005}. To speed up the process, the interpolation coefficients were precalculated and stored in a look-up table. In the next step the interpolated data were normalized such that the $L^2$ norm equaled 100. This allows for choosing the reconstruction parameters independent from the data acquisition parameters, which minimizes the operator interference and also maintains the quality of the results.

For AME reconstruction of experimental datasets, we found empirically that it is sufficient to choose $K = \{-2, -1, 0, 1, 2\}$ to exploit the complementary information from 5 successive acquisitions with interleaved radial encodings. Numerically, $K = \{-3, -2, -1, 0, 1, 2, 3\}$ gives similar results for both simulated and real data (not shown) but increases computational complexity. The preliminary images were reconstructed by NLINV, with almost identical parameters for regularization and penalization of coil sensitivity profiles. The same initialization was used with the spin density set to ones and coil sensitivities to zeros for the first frame. Later both were replaced by the reconstruction results from the previous frame. For motion estimation, the spatial deformation of the object was implemented with a bicubic interpolation, and parameters of the model were set to $\lambda = 0.02$ and $\mu = 1.0$.

For comparison, the same data was also reconstructed by the standard NLINV method \cite{Uecker2010a,Zhang2010a} which was implemented in the same software environment. In addition, a temporal median filter was applied to the images to reduce residual streaking artifacts as proposed in \cite{Uecker2010,Uecker2010a,Zhang2010a}. It was implemented as a post-processing step with a window width covering 5 neighboring images (NLINV-MED).

\section{Results}

\subsection{Motion Estimation}
\begin{figure}[!t]
\centering
\includegraphics[width=0.8\textwidth]{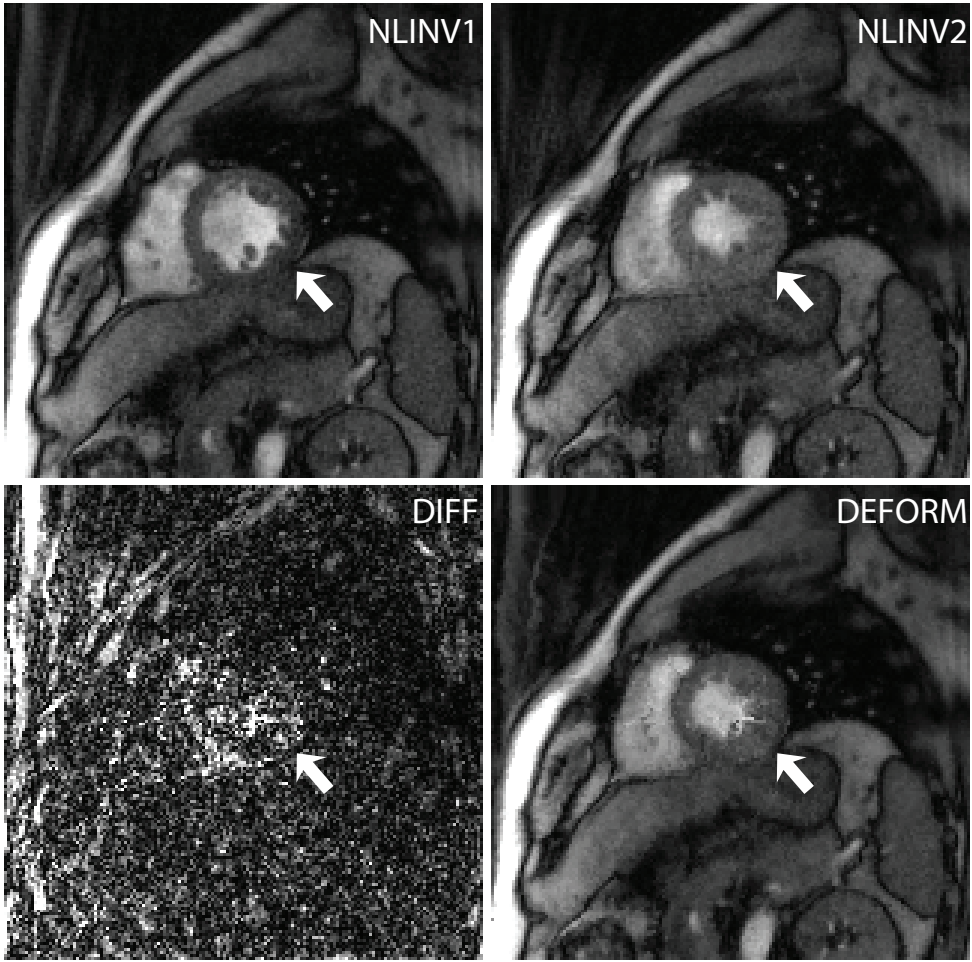}
\caption{Principle performance of the motion estimation algorithm (short-axis cardiac MRI). For two frames at end diastole (NLINV1) and end systole (NLINV2) motion fields were estimated based on their NLINV reconstructions and then applied to NLINV1 as DEFORM. DIFF represents the difference between NLINV2 and DEFORM. }
\label{motionest}
\end{figure}
The principle of the proposed motion estimation is demonstrated in Figure \ref{motionest} using data for the human heart. Two frames at end diastole (NLINV1) and end systole (NLINV2) were selected to depict distinct differences due to myocardial contraction in preliminary NLINV reconstructions (arrows). The calculated deformation of NLINV1 by incorporating the estimated motion is shown as DEFORM. It clearly identifies the contraction of the myocardium, whereas the streaking artifacts at the top-left corner of the image remain similar as in NLINV1. The example demonstrates that the information of the moving object has correctly been captured by the motion estimation, while the image artifacts, which may also change from time to time, are appropriately excluded. This can also be visualized in the difference between DEFORM and NLINV2 (DIFF), where the artifacts are dominant and the structure of the heart is less visible. In the next subsection, we show how the artifacts will be removed, rather than enhanced, by aggregating the estimated motions in AME (cf. Figure \ref{workflow}).

\subsection{AME in Action}
\subsubsection{Motion Phantom}
\begin{figure}[!t]
\centering
\includegraphics[width=0.8\textwidth]{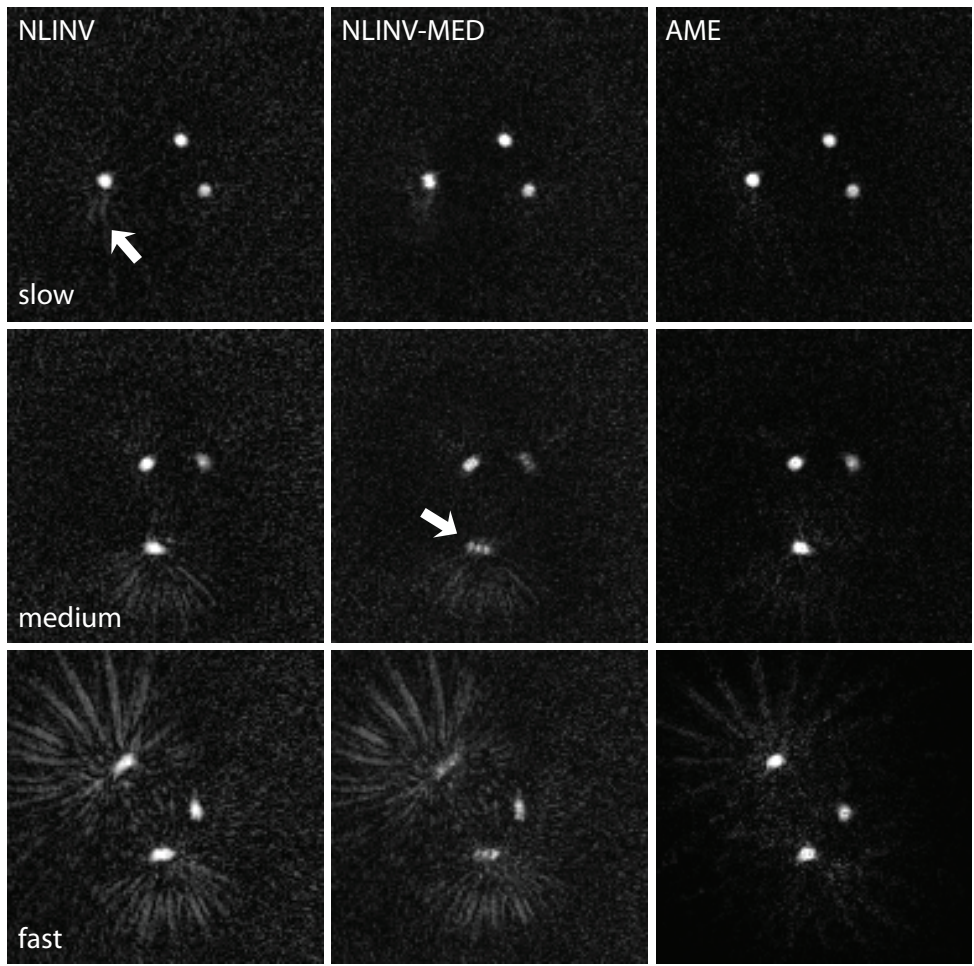}
\caption{Reconstructions from a phantom (three tubes of water) moving with different constant angular speed for (left) NLINV, (middle) NLINV-MED, and (right) AME. The rotating frequencies of the moving tubes are (slow) 0.5 Hz, (medium) 1.0 Hz and (fast) 1.5 Hz, respectively. }
\label{movphan}
\end{figure}
Figure \ref{movphan} compares reconstructions for a phantom moving at different speeds that were obtained by NLINV, NLINV-MED, and the proposed AME method, respectively. For the lowest velocity, all three methods produce acceptable results, although the latter two surpass NLINV in reducing streaking artifacts (arrow in top row). At moderate velocity, NLINV suffers from stronger artifacts due to faster motions, while NLINV-MED even distorts the structure of the fastest moving outermost tube (arrow in middle row). The stretched shape of the circular tube is a typical effect from the temporal median filter. In contrast, the AME reconstruction offers a proper image with almost no motion or streaking artifacts. Finally, for the highest velocity, both NLINV and NLINV-MED result in severely deformed shapes for almost all tubes as well as pronounced streaking artifacts. Again, the AME method shows best results with only very mild and barely visible artifacts. Furthermore, the signal-to-noise ratio (SNR) of the AME reconstruction is higher in all cases compared with the two other methods.

\subsubsection{Human Tongue Movements}
\begin{figure}[!t]
\centering
\includegraphics[width=0.6\textwidth]{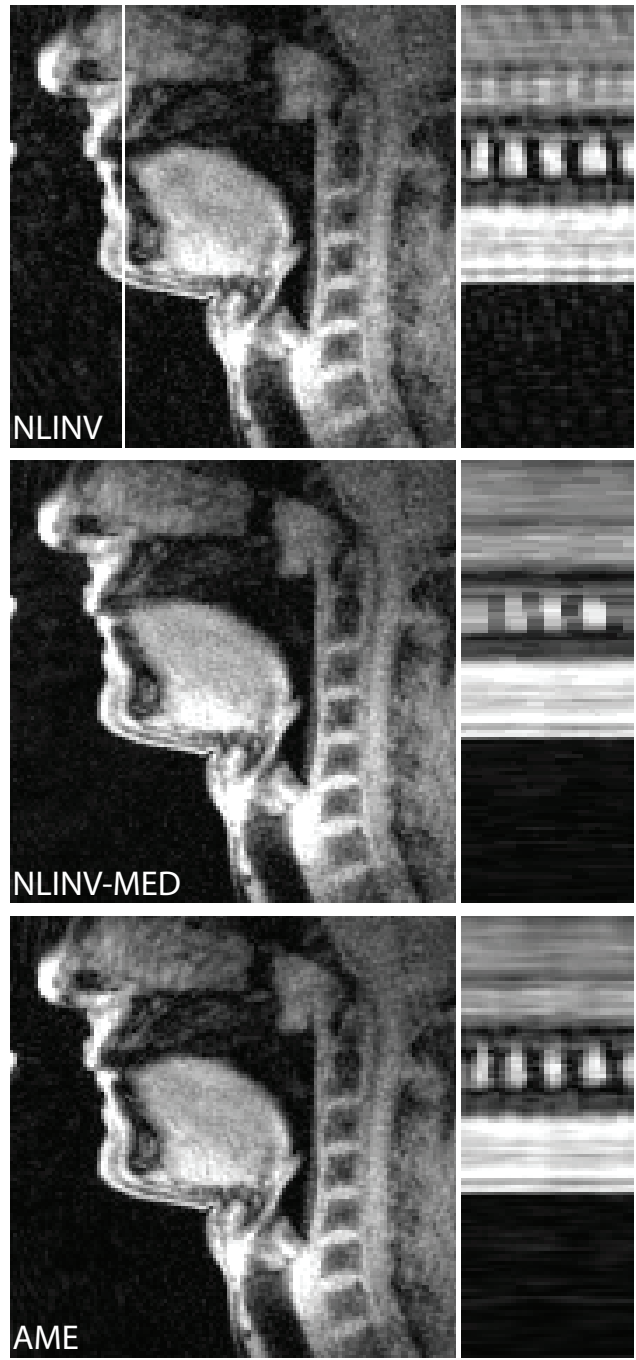}
\caption{Reconstructions of tongue movements in a mid-sagittal plane at a temporal resolution of 33 ms for NLINV, NLINV-MED, and the proposed AME. The right column presents spatiotemporal profiles (vertical line in upper left image) of the MRI signal intensity for a selected 1.0 s period.}
\label{tongue_motion}
\end{figure}
Figure \ref{tongue_motion} demonstrates tongue movements during playing the mouthpiece of a brass instrument. For this particular task the tongue tip of the subject had to rapidly move forward and backward touching the upper teeth ridge. To better demonstrate the temporal evolution of the motion, a reference line is placed at the tongue tip to derive corresponding 2D spatiotemporal intensity profiles. Thus, the flickering of the residual artifacts at every 5-th frame is clearly visualized for NLINV. For NLINV-MED the residual artifacts are effectively removed at the expense of blurring the tongue movements by the temporal median filter. On the contrary, the proposed AME method preserves the sharp intensity changes associated with the rapid tongue movements even better than in the original NLINV reconstruction, while at the same time successfully minimizing residual streaking artifacts.

\section{Discussion}
In comparison to NLINV reconstructions with and without temporal median filter, the proposed AME reconstruction for real-time MRI with pronounced radial undersampling yields serial images with improved temporal acuity and less residual artifacts. The new approach emerges as an expansion of the previously introduced NLINV reconstruction with an aggregated motion estimation which estimates respective movements from multiple consecutive data sets with complementary spatial encodings. The additional information is incorporated into the data consistency term of the nonlinear inverse problem for a simultaneous determination of spin density and coil sensitivities.

Extending other approaches for motion estimation in MRI, the present work is not limited to affine or rigid motions. Moreover, the combination of AME with nonlinear reconstruction permits an arbitrary choice of $K$ which defines the set of frames used for reconstructing the actual frame. Extensive experimental studies (not shown) demonstrate that a choice of $|K|$ smaller than the number of frames with complementary spatial encodings fails to remove the temporally flickering artifacts. On the other hand, choosing $|K|$ greater than the number of differently encoded frames does not further improve the reconstruction but yields comparable image quality. Because the computational complexity increases as $|K|$ increases, this behavior explains our choice of $K = \{-2,-1,0,1,2\}$ in all experiments. As a stopping criterion of the iterations in NLINV, NLINV-MED and AME, the well-known Morozov's discrepancy principle was initially considered, but it forbids a unique choice of the threshold value for every frame because the energy of the signal slightly changes with time even for normalized $k$-space data. In our applications, we have chosen a fixed number of iterations (i.e., Newton steps) for each method, respectively, which gives satisfactory results. A data-driven choice certainly appears to be more sensible and might be considered in future research.

A limitation of the AME method, which may deteriorate its performance, stems from errors in the preceding NLINV reconstructions that may lead to an unfaithful motion estimation. This can be seen from the slightly blurred temporal profile in Figure \ref{tongue_motion}. A natural way to overcome this problem would be to run the whole AME procedure at least twice using previous AME (rather than NLINV) reconstructions for more accurate motion estimations, though at the expense of further increasing the computational demand. In fact, at this time the high computational cost, which is about 20 times that of a comparable NLINV implementation on a laptop with MATLAB, is currently the major obstacle for more extended practical applications. However, because the computations for each receiver coil and of different frames are independent, AME is highly adaptable to parallel computing. Apart from interpolation, the involved calculations are simplified to point-wise operations, fast Fourier transform, and scalar products. As is shown in Section II-D, the interpolation for non-Cartesian data may be separated from the iterative optimization, through a convolution with the point-spread function. These features further ensure a possible speed-up by an implementation on graphical processing units.

\section{Conclusion}
This work introduces a new reconstruction method for real-time MRI that offers improved temporal fidelity for visualizing rapid dynamic changes. Preliminary results for an experimental phantom and \textit{in vivo} human data demonstrate the practical performance and improved quality which is based on the incorporation of estimated object motions into the nonlinear inverse reconstruction process. Future improvements are expected by exploiting new regularization methods and by accelerating the computational speed.

\appendix

\section*{Acknowledgment}
The authors thank Sebastian Schaetz and Aaron Niebergall for the design and construction of the motion phantom. H. Li  thanks for financial support by the China Scholarship Council. S. Zhang thanks for financial support by the DZHK (German Centre for Cardiovascular Research) and BMBF (German Ministry of Education and Research). A. Munk gratefully acknowledges the support of DFG FOR 916.


\end{document}